\documentclass[12pt,a4paper]{amsart}

\usepackage{euscript,amsfonts,amssymb,amsmath,amscd,amsthm,enumerate,hyperref}
\usepackage{comment}
\usepackage{mathrsfs}
\usepackage{graphicx}
\usepackage{xcolor}
\usepackage[textsize=tiny]{todonotes}
\usepackage[margin=1.1in]{geometry}
\usepackage{bbm}
\usepackage{enumitem}
\usepackage{upgreek}
\usepackage{enumerate}
\usepackage[giveninits=true,style=alphabetic,backend=biber,maxnames=9,maxalphanames=9,doi=false,isbn=false,url=false,eprint=false]{biblatex}

\renewbibmacro{in:}{}
\AtEveryBibitem{%
 \clearlist{language}%
 \clearfield{note}%
}

\newcommand{\lb}{\left}
\newcommand{\rb}{\right}

\newcommand{\R}{\mathbb{R}}

\newcommand{\dd}{\mathrm{d}}

\newcommand{\veps}{\varepsilon}
\newcommand{\sgn}{\mathrm{sgn}}

\newcommand{\wt}{\widetilde}
\newcommand{\wh}{\widehat}

\newcommand{\calF}{\mathcal{F}}
\newcommand{\calP}{\mathcal{P}}
\newcommand{\calG}{\mathcal{G}}
\newcommand{\calR}{\mathcal{R}}

\newcommand{\scrX}{\mathscr{X}}
\newcommand{\scrY}{\mathscr{Y}}

\newcommand{\bP}{\mathbb{P}}

\newcommand{\E}{\mathbb{E}}

\newcommand{\ind}{\mathbf{1}}

\usepackage{amsthm}

\theoremstyle{plain}
\newtheorem{thm}{Theorem}[section]
\newtheorem*{thm*}{Theorem}

\newtheorem{lem}[thm]{Lemma}
\newtheorem{ass}[thm]{Assumption}

\theoremstyle{definition}
\newtheorem{defn}[thm]{Definition}

\newtheorem{rmk}[thm]{Remark}

\title[Particle Systems with Interaction]{A singular Two-Phase Stefan problem and \\ particles Interacting Through Their Hitting Times}

\author{Graeme Baker}
\address{Program in Applied \& Computational Mathematics, Princeton University, Princeton, NJ 08544, USA}
\email{graemeb@math.princeton.edu}

\author{Mykhaylo Shkolnikov}
\address{ORFE Department, Bendheim Center for Finance, and Program in Applied \& Computational Mathematics, Princeton University, Princeton, NJ 08544, USA}
\email{mshkolni@gmail.com}
\thanks{M.~Shkolnikov is partially supported by the NSF grant DMS-2108680.}

\begin{document}

\begin{abstract}
We consider a probabilistic formulation of a singular \textit{two-phase} Stefan problem in one space dimension, which amounts to a coupled system of two McKean-Vlasov stochastic differential equations.
~In the financial context of systemic risk, this system models two competing regions with a large number of interconnected banks or firms at risk of default.
~Our main result shows the existence of a solution whose discontinuities obey the natural physicality condition for the problem at hand.~Thus, this work extends the recent series of existence results for singular \textit{one-phase} Stefan problems in one space dimension that can be found in \cite{delarue_particle_2015}, \cite{nadtochiy_mean_2019}, \cite{hambly_mckean--vlasov_2018}, \cite{cuchiero_propagation_2020}.~As therein, our existence result is obtained via a large system limit of a finite particle system approximation in the Skorokhod M1 topology. But, unlike for the previously studied one-phase case, the free boundary herein is not monotone, so that the large system limit is obtained by a novel argument. 
\end{abstract}

\maketitle

\section{Introduction}

We are concerned with the following problem of McKean-Vlasov type. Given random variables $X_{0-}\ge 0$, $Y_{0-}\le 0$ and parameters $\alpha,\beta>0$ one seeks a right-continuous function with left limits $\Lambda\!:[0,\infty)\to\R$ that solves the system of equations
\begin{equation}\label{limitprob}
\begin{split}
&X_t=X_{0-}+B_t-\Lambda_t,\\
&Y_t=Y_{0-}+W_t-\Lambda_t,\\
&\Lambda_t=\alpha \bP(\sigma\le t)-\beta \bP(\tau\le t),
\end{split}
\end{equation}
where $B$, $W$ are standard Brownian motions independent of $X_{0-}$, $Y_{0-}$, respectively, $\sigma:=\inf\{t\ge0\!:X_t\le 0\}$, and $\tau:=\inf\{t\ge0\!:Y_t\ge 0\}$.

\medskip

Problem \eqref{limitprob} with $\beta=0$ and its variants, in which the process $Y$ can be omitted, have received much attention recently. More specifically, they have been shown to arise in the large scale limit of integrate-and-fire models from computational neuroscience (see \cite{delarue_global_2015}, \cite{delarue_particle_2015}), of models for borrowing and lending between financial institutions (see \cite{nadtochiy_particle_2019}), and of random growth models from mathematical physics (see \cite{dembo_criticality_2019}). In addition, if $X_{0-}$ admits a density $f$ in the Sobolev space $W^1_2([0,\infty))$ with $f(0)=0$ and the derivative $\Lambda'$ exists as a function in $L^2([0,T])$ for some $T>0$, then for each $t\in[0,T]$ the random variable $(X_t+\Lambda_t)\,\mathbf{1}_{\{\sigma>t\}}$ possesses a density $u(t,\cdot)$ on $(\Lambda_t,\infty)$, and the pair $(\Lambda,u)$ is a classical solution of the one-phase supercooled Stefan problem for the heat equation on $[0,T]$ (cf.~\cite[proof of Proposition~4.2(b)]{nadtochiy_particle_2019}, and see \cite[Remark 1.2]{delarue_global_2022} for a stronger version of this statement).~The latter is the prototypical model for the solidification of a supercooled liquid and has been studied extensively since the 1970's (see, e.g., the survey \cite{visintin_introduction_1998}). 

\medskip

The arguments in \cite[proof of Proposition~4.2(b)]{nadtochiy_particle_2019} also apply to problem \eqref{limitprob} in its full generality. To wit, if $X_{0-}$, $Y_{0-}$ admit densities $f$, $g$ in $W^1_2([0,\infty))$, $W^1_2((-\infty,0])$, respectively, with $f(0)\!=\!g(0)\!=\!0$ and $\Lambda'$ exists as a function in 
$L^2([0,T])$ for some $T\!>\!0$, then for each $t\in[0,T]$
the random variables $(X_t+\Lambda_t)\,\mathbf{1}_{\{\sigma>t\}}$, $(Y_t+\Lambda_t)\,\mathbf{1}_{\{\tau>t\}}$ possess densities $u(t,\cdot)$, $v(t,\cdot)$ on $(\Lambda_t,\infty)$, $(-\infty,\Lambda_t)$, respectively, and the triple $(\Lambda,u,v)$ is a classical solution of the \textit{two-phase} Stefan problem 
\begin{equation}\label{2phasePDE}
\begin{split}
& \partial_t u = \frac{1}{2}\partial_{xx} u\;\;\;\text{on}\;\;\;\{(t,x)\in[0,T]\times\R:\,x\ge\Lambda_t\}, \\
& \partial_t v = \frac{1}{2}\partial_{xx} v\;\;\;\text{on}\;\;\;\{(t,x)\in[0,T]\times\R:\,x\le\Lambda_t\}, \\
& \Lambda'_t = \frac{\alpha}{2}\partial_x u(t,\Lambda_t)+\frac{\beta}{2}\partial_x v(t,\Lambda_t),\quad t\ge0, \\
& u(0,x)=f(x),\;\;\;x\ge0\quad\text{and}\quad u(t,\Lambda_t)=0,\;\;\; t\ge0, \\
& v(0,x)=g(x),\;\;\;x\le0\quad\text{and}\quad v(t,\Lambda_t)=0,\;\;\; t\ge0. 
\end{split}
\end{equation}
\smallskip

In the language of phase change problems, \eqref{2phasePDE} describes the process of freezing and melting over time when a liquid with the initial temperature distribution of $-f$ (relative to its equilibrium freezing point) is placed next to the corresponding solid with the initial temperature distribution of $g$ (again, relative to the equilibrium freezing point).
~To clarify terms, the liquid is said to be supercooled since it has been cooled below the freezing point, and similarly the solid is superheated as its temperature has been raised above the melting point.
~In this framework, $\Lambda_t$ gives the location of the interface between the two phases and $-u(t,\cdot)$, $v(t,\cdot)$ provide the temperature distribution in the two phases at time $t$.~The physically more relevant situation where a supercooled liquid is put in contact with the corresponding non-superheated solid was studied in \cite{gotz_two-phase_1995}.~Therein, the authors show the existence of a weak solution and pin down, in a case where the interface has at most one jump, the size of that jump.~In the notation of~\eqref{2phasePDE}, their setting corresponds to taking $\alpha\ge0\ge\beta$ which leads to a non-decreasing interface $\Lambda$.~We also refer to \cite[(1.12)--(1.14)]{meirmanov_1994} where the same sign convention arises from a Stefan problem with surface tension in $\R^3$ under a radial symmetry assumption.~While our techniques apply also for $\alpha\!\ge\!0\!\ge\!\beta$ (see Remark \ref{rmk:alphabeta}), we focus on the case~$\alpha,\beta>0$ which results in a \textit{non-monotone} interface $\Lambda$, making the mathematical analysis more challenging.~In this case, \eqref{2phasePDE} can be viewed as a linear one-dimensional prototype for the Stefan problems conjectured to arise in the sharp interface limit of the Giacomin-Lebowitz phase segregation model (cf.~\cite[Statement~1]{GL}).

\medskip

Beyond mathematical physics, we hope that the solutions of the two-phase problems \eqref{limitprob} and \eqref{2phasePDE} constructed herein will prove useful for applications in other areas, such as neuroscience (e.g.,~to integrate-and-fire models featuring both excitatory and inhibitory neurons) and mathematical finance. With regards to the latter, in Section \ref{subsec:app} we build on the systemic risk model of \cite{nadtochiy_particle_2019} by introducing a model for two competing regions with interconnected financial firms at risk of default. In this context, $\Lambda_t$ is determined by the proportion of firms from each region which have defaulted by time $t$. 
As the number of firms grows large, our main result shows that the representative asset values will satisfy equation \eqref{limitprob}. In Section \ref{sec:conc}, we discuss the systemic risk implications of this model and propose possible extensions.

\medskip

Returning to the probabilistic formulation \eqref{limitprob}, we note first that its solutions $\Lambda$ fail to be continuous for generic $X_{0-}\ge0$, $Y_{0-}\le0$, $\alpha>0$, $\beta>0$.~For example, adapting \cite[proof of Theorem 1.1]{hambly_mckean--vlasov_2018} to the setting of \eqref{limitprob} we see that $\Lambda$ must be discontinuous whenever $\alpha\E[X_{0-}]-\beta\E[Y_{0-}]<(\alpha-\beta)^2/2$. Indeed, assuming the opposite we infer from $X_{t\wedge\sigma}=X_{0-}+B_{t\wedge\sigma}-\Lambda_{t\wedge\sigma}\ge0$ and $Y_{t\wedge\tau}=Y_{0-}+W_{t\wedge\tau}-\Lambda_{t\wedge\tau}\le 0$ that
\begin{eqnarray}
&& \E[X_{0-}]\ge \E[\Lambda_{t\wedge\sigma}]\underset{t\to\infty}{\longrightarrow} \E[\Lambda_\sigma]
=\int_0^\infty \Lambda_s\,\mathrm{d}\bP(\sigma\le s), \label{ineq1} \\
&& \E[Y_{0-}]\le \E[\Lambda_{t\wedge\tau}]\underset{t\to\infty}{\longrightarrow} \E[\Lambda_\tau]
=\int_0^\infty \Lambda_s\,\mathrm{d}\bP(\tau\le s). \label{ineq2}
\end{eqnarray}
Combining \eqref{ineq1}, \eqref{ineq2}, the third equation in \eqref{limitprob}, and $\Lambda_0=0$ we deduce that
\begin{equation}
\alpha \E[X_{0-}] - \beta \E[Y_{0-}] \ge \int_0^\infty \Lambda_s\,\mathrm{d}\Lambda_s = \frac{\Lambda_\infty^2}{2}. \label{ineq3}
\end{equation}
Since $\Lambda$ takes values in $[-\beta,\alpha]$, we have $\Lambda_\infty=\alpha\bP(\sigma<\infty)-\beta\bP(\tau<\infty)=\alpha-\beta$, so \eqref{ineq3} poses the desired contradiction to the assumed $\alpha\E[X_{0-}]-\beta\E[Y_{0-}]<(\alpha-\beta)^2/2$. 

\medskip

We prove in Section \ref{sec:limit} (see the proof of Theorem \ref{thm:exist}(iv), lower bound) that the jump sizes $|\Lambda_t-\Lambda_{t-}|:=|\Lambda_t-\lim_{s\uparrow t} \Lambda_s|$ in any solution $\Lambda$ of \eqref{limitprob} must satisfy the inequality 
\begin{equation}\label{jumplowerbound}
|\Lambda_t-\Lambda_{t-}|\ge
\begin{cases}
\inf\big\{x>0:\,\alpha \bP(\sigma\ge t,X_{t-}\in [0,x])<x\big\}&\text{if}\quad\Lambda_t-\Lambda_{t-}>0,\\
\inf\big\{y>0:\,\beta \bP(\tau\ge t,Y_{t-}\in [-y,0])<y\big\}&\text{if}\quad\Lambda_t-\Lambda_{t-}<0.
\end{cases}
\end{equation}
This motivates to search for solutions of \eqref{limitprob} for which \eqref{jumplowerbound} holds with equality.~In analogy with the corresponding terminology for the one-phase case $\beta=0$ we refer to such solutions as \textit{physical}. 

\begin{defn}\label{def:phys}
A solution $\Lambda$ of \eqref{limitprob} is said to be \emph{physical} if for any $t$ with $\Lambda_t\neq\Lambda_{t-}$, 
\begin{equation}\label{eq:physcond}
|\Lambda_t-\Lambda_{t-}|=
\begin{cases}
	\inf\big\{x>0:\,\alpha \bP(\sigma\ge t,X_{t-}\in [0,x])<x\big\}&\text{if}\quad\Lambda_t-\Lambda_{t-}>0,\\
	\inf\big\{y>0:\,\beta \bP(\tau\ge t,Y_{t-}\in [-y,0])<y\big\}&\text{if}\quad\Lambda_t-\Lambda_{t-}<0.
\end{cases}
\end{equation}
\end{defn}

\begin{rmk}
An interesting feature of the solutions to \eqref{limitprob} we construct is that a positive jump in $\Lambda$ at any given time rules out the possibility of a jump in $\beta \bP(\tau\le t)$ at that same time (and similarly with negative jumps in $\Lambda$ ruling out simultaneous jumps in $\alpha \bP(\sigma\le t)$). Therefore, it is natural to choose the jump sizes of physical solutions in each phase according to the same physicality condition as in the one-phase case.
\end{rmk}

We are now ready to state our main result. 

\begin{thm*}
\!\! Let $\alpha,\beta\!>\!0$.~If $\,\E[X_{0-}]\!<\!\infty$, $\E[Y_{0-}]\!>\!-\infty$ and $\bP(X_{0-}\!=\!0)\!=\!0\!=\!\bP(Y_{0-}=0)$, then there exists a physical solution $\Lambda$ of \eqref{limitprob}.
\end{thm*}

\begin{rmk}
The integrability assumption on the initial condition in the theorem matches the currently best known assumption under which existence of a physical solution has been established in the one-phase case $\beta=0$ (see \cite[Theorem~6.5]{cuchiero_propagation_2020}). The additional assumption that $0$ is not an atom is natural, since there does not seem to be a canonical way to attribute such an atom to one or to both of the phases. 
\end{rmk}

Our proof of the theorem is based on a finite particle system approximation to~\eqref{limitprob}, which extends the ones invented for the one-phase case $\beta=0$ (see \cite{delarue_particle_2015}, \cite{nadtochiy_particle_2019}, \cite{hambly_mckean--vlasov_2018}, \cite{cuchiero_propagation_2020}).~In the finite particle systems herein, when a particle from a given phase hits the boundary, the particles in the same phase shift towards the boundary and the particles in the other phase shift away from it.~When many particles from the same phase are close to the boundary and one hits, this may trigger a cascade where a macroscopic number of particles hit the boundary instantaneously and all particles from both phases jump a proportionally large amount.~Similar cascade events correspond to important phenomena in the applications mentioned above, such as the synchronization of neurons, the systemic default of a large fraction of financial institutions, and the instantaneous freezing of a supercooled liquid.

\medskip

As the number of particles goes to infinity, the cascades do not disappear, but rather manifest themselves as jumps in the movement of a representative particle. For the two-phase model herein, the cascades originating from the two phases compete and interact with each other.~This feature is not present in the one-phase case and necessitates new techniques to show the existence of physical large system limits.~The McKean-Vlasov problem~\eqref{limitprob} describing these limits may admit multiple solutions, corresponding to different choices of the jump sizes.~However, for the finite particle systems there is a single natural way to resolve cascades that respects causality, and this leads to a characterization of the jump sizes in terms of the empirical measure of the particles right before the jumps.~This further motivates the definition of physical solutions in terms of the laws of the representative particles $X$ and $Y$ (Definition \ref{def:phys}).

\medskip

We note the possibility that both $\inf\{ x>0\!: \alpha \bP(\sigma\ge t,X_{t-}\in [0,x])<x\}>0$ and $\inf\{y>0\!: \beta \bP(\tau\ge t,Y_{t-}\in [-y,0])<y\}>0$ at some $t$.~For instance, this can happen at time $t=0$ for suitable choices of the distributions of $X_{0-}$ and $Y_{0-}$, or at any time $t>0$ when $\alpha=\beta$ and $Y_{0-}\stackrel{d}{=}-X_{0-}$.~When this occurs, there are multiple physical solutions to \eqref{limitprob} because our main theorem applies on $[t,\infty)$ once a choice of $\Lambda_t$ has been made. It remains an open problem to find assumptions on the initial condition which guarantee the uniqueness of the physical solution to the two-phase problem.~Uniqueness for the one-phase problem (see \cite{delarue_global_2022}) relies on a comparison principle which is not immediately applicable here.
As in the setting of \cite{delarue_global_2022}, the laws of the representative particles $X$ and $Y$ in \eqref{limitprob} give a global notion of solution to the corresponding free boundary problem \eqref{2phasePDE} and allow to continue past the times of blow-up in $\Lambda'$.~If, in addition, the physical solution turns out to be unique, this would yield the well-posedness of the underlying free boundary problem, despite the possible presence of such blow-ups. 

\medskip

In contrast to the one-phase problem, the free boundary $\Lambda$ in \eqref{limitprob} is not necessarily monotone, as it is the difference between the competing terms $\alpha \bP(\sigma\!\le\! t)$ and $\beta \bP(\tau\!\le\! t)$. Nonetheless, with new arguments including \cite[proof idea of Proposition 3.4(i)]{antunovic_isolated_2011}, we show the existence of a physical solution to \eqref{limitprob} by taking the large system limit of a suitable finite particle system, which is introduced in Section \ref{sec:part}. In Section \ref{sec:tight}, we prove the tightness of the laws of physical solutions to such systems as the number of particles goes to infinity.~Then, in Section \ref{sec:limit}, we check that the limiting law of a representative particle gives a physical solution to \eqref{limitprob}. 

\medskip

Throughout the paper, $\varphi$ and $\Phi$ denote the standard Gaussian probability density function (pdf) and cumulative distribution function (cdf), respectively. We use the term c\`{a}dl\`{a}g to refer to functions that are right-continuous with left limits. In addition, we use the notation $h(t-)$ for the left limit $\lim_{s\uparrow t} h(s)$ whenever we are working with a c\`adl\`ag function $h$. 

\section{Interacting Particle System Approximation}\label{sec:part}

We obtain the existence of physical solutions to \eqref{limitprob} by considering the large system limit of the interacting particle system introduced in this section. We are going to see that the physicality of the limit is inherited from the way that we resolve the cascades (i.e., the jumps in the interaction term) for the particle system at size $N$. Only after constructing the solutions in this section do we give the equation for the common drift $\Lambda^N$ in terms of hitting times. The reason for this is that the equations admit non-physical solutions, where $\Lambda^{N}$ jumps, causing particles to hit, justifying the jump in $\Lambda^{N}$; however, this situation does not respect causality. We refer the reader to the related discussion of physical and non-physical solutions in \cite[Subsection 3.1]{delarue_particle_2015}.

\medskip

For $N\ge 1$, let the coupled particle systems $(X^{N,1},\dots,X^{N,N})$ and $(Y^{N,1},\dots,Y^{N,N})$ be given by
\begin{equation}\label{eq:partsys}
\begin{split}
X_t^{N,i}&=X_{0-}^{N,i}+B_t^{i}-\Lambda_t^{N},\\
Y_t^{N,i}&=Y_{0-}^{N,i}+W_t^{i}-\Lambda_t^{N},
\end{split}
\end{equation}
where $B^{1},\dots,B^{N}$ and $W^{1},\dots,W^{N}$ are independent standard Brownian motions and the initial conditions $X_{0-}^{N,1},\dots,X_{0-}^{N,N}$ and $Y_{0-}^{N,1},\dots, Y_{0-}^{N,N}$ are drawn independently from the laws $\mu_{0-}$ and $\nu_{0-}$ of $X_{0-}$ and $Y_{0-}$, respectively. Note that we keep track of the dependence on the size of the systems $N$ in the superscript of $X^{N,i}$ and $Y^{N,i}$; however, we may use the same independent Brownian motions for all $N$. For the initial data, we make the same assumption as in the main theorem.
\begin{ass}\label{ass1}
Let both $\mu_{0-}\in\calP([0,\infty))$ and $\nu_{0-}\in\calP((-\infty,0])$ have a finite first moment and no atom at $0$.
\end{ass}

We stipulate that $\Lambda^{N}$ is zero up until the first cascade time
\begin{align}
T_1:=\inf\big\{t\ge0:\,X_{t-}^{N,i}\le 0\text{ or }Y_{t-}^{N,i}\ge 0\text{ for some }i\in\{1,\dots,N\}\big\}.
\end{align}
Following \cite[Subsection 3.1]{delarue_particle_2015}, we introduce a second time axis, the \emph{cascade time axis at} $T_1$. At the start of the cascade time axis we consider the sets
\begin{equation}
\begin{split}
\Sigma_{1}^0=\big\{i\in\{1,\dots,N\}:\,X_{T_1-}^{N,i}\le 0\big\}\quad\text{and}\quad
\Gamma_{1}^0=\big\{i\in\{1,\dots,N\}:\,Y_{T_1-}^{N,i}\ge 0\big\}.
\end{split}
\end{equation}
These are the particles which reach $0$ ``before $\Lambda^{N}$ jumps'', activating the cascade. As $\mu_{0-}$ and $\nu_{0-}$ do not put mass at the origin, with probability one, only one of $\Sigma_1^{0}$ and $\Gamma_1^{0}$ will be non-empty and it will be a singleton. As we wish to work on the whole probability space $(\Omega,\calF, \bP)$, we disregard this and stipulate what happens when multiple particles reach the boundary simultaneously. Considering this case is necessary in the extended problem below (see Remark \ref{rmk:comparison}) where each particle's motion is driven by an arbitrary c\`adl\`ag process.

\medskip

We define
\begin{equation}
\begin{split}
\Sigma_{1}^1&=\Big\{i\in\{1,\dots,N\}\setminus\Sigma_1^{0}:\,X_{T_1-}^{N,i}\le \frac{\alpha}{N}|\Sigma_{1}^0|-\frac{\beta}{N}|\Gamma_1^{0}|\Big\}\quad\text{and}\\
\Gamma_{1}^1&=\Big\{i\in\{1,\dots,N\}\setminus\Gamma_1^{0}:\,Y_{T_1-}^{N,i}\ge \frac{\alpha}{N}|\Sigma_{1}^0|-\frac{\beta}{N}|\Gamma_{1}^0|\Big\}.
\end{split}
\end{equation}
Moving by one unit of cascade time, the sets $\Sigma_{1}^1$ and $\Gamma_{1}^1$ represent further particles that cross zero due to the cascade initiated by the particles in $\Sigma_{1}^0$ and $\Gamma_{1}^0$. We note that $\Sigma_{1}^1$ will be empty if $\frac{\alpha}{N}|\Sigma_{1}^0|-\frac{\beta}{N}|\Gamma_1^{0}|$ is non-positive and $\Gamma_{1}^1$ will be empty if $\frac{\alpha}{N}|\Sigma_{1}^0|-\frac{\beta}{N}|\Gamma_1^{0}|$ is non-negative. Similarly, for $k=2,\dots, N-1$, we set
\begin{equation}
\begin{split}
\Sigma_{1}^k&=\Big\{i\in\{1,\dots,N\}\setminus\cup_{j=0}^{k-1}\,\Sigma_1^{j}:\,X_{T_1-}^{N,i}\le \frac{\alpha}{N}\big|\!\cup_{j=0}^{k-1} \Sigma_1^{j}\big|
-\frac{\beta}{N}\big|\!\cup_{j=0}^{k-1} \Gamma_1^{j}\big|\Big\}\quad
\text{and} \\
\Gamma_{1}^k&=\Big\{i\in\{1,\dots,N\}\setminus\cup_{j=0}^{k-1}\,\Gamma_1^{j}:\,Y_{T_1-}^{N,j}\ge \frac{\alpha}{N}\big|\!\cup_{j=0}^{k-1}\Sigma_1^{j}\big|
-\frac{\beta}{N}\big|\!\cup_{j=0}^{k-1}\Gamma_1^{j}\big|\Big\}.
\end{split}
\end{equation}
If $\Sigma_{1}^k$ or $\Gamma_1^k$ is empty for some $k$ then it will be empty for all $\ell>k$. We may stop at $\Sigma_{1}^{N-1}$ (and $\Gamma_{1}^{N-1}$) since in the most extreme scenario each of $\Sigma_{1}^{0},\dots,\Sigma_1^{N-1}$ (or $\Gamma_{1}^{0},\dots,\Gamma_1^{N-1}$) has exactly one element.

\medskip

We set $\Lambda_{T_1}^N=\frac{\alpha}{N}\big|\!\cup_{j=0}^{N-1}\Sigma_1^{j}\big|-\frac{\beta}{N}\big|\!\cup_{j=0}^{N-1}\Gamma_1^{j}\big|$ and stipulate that $\Lambda^N$ is constant until the second cascade time
\begin{equation}\label{eq:T2}
\begin{split}
T_2:=\inf\big\{t>T_1:\;\,&X_{t-}^{N,i}\le 0\text{ for some }i\in\{1,\dots,N\}\setminus\cup_{j=0}^{N-1}\Sigma_1^{j}\\
&\text{or } Y_{t-}^{N,i}\ge 0\text{ for some }i\in\{1,\dots,N\}\setminus\cup_{j=0}^{N-1}\Gamma_1^{j}\big\}.
\end{split}
\end{equation}
In the definition of $T_2$ we are careful not to include particles that have already crossed zero. Proceeding along the cascade time axis at $T_2$, the sets $\Sigma_2^{0},\dots,\Sigma_2^{N-1}$ and $\Gamma_2^{0},\dots,\Gamma_2^{N-1}$ are defined analogously to $\Sigma_1^{0},\dots,\Sigma_1^{N-1}$ and $\Gamma_1^{0},\dots,\Gamma_1^{N-1}$, where we again do not include particles that have already crossed zero. We define $\Lambda_{T_2}^N=\Lambda_{T_1}^N+\frac{\alpha}{N}\big|\!\cup_{j=0}^{N-1}\Sigma_2^{j}\big|-\frac{\beta}{N}\big|\!\cup_{j=0}^{N-1}\Gamma_2^{j}\big|$ and continue the process \emph{mutatis mutandis}.

\medskip

With $\Lambda^N$ defined we can introduce the stopping times
\begin{align}\label{eq:sigmadef}
\sigma^{N,i}:=\inf\{t\ge0:\; X_{t}^{N,i}\le 0\}\quad\text{and}\quad
\tau^{N,i}:=\inf\{t\ge0:\;Y_{t}^{N,i}\ge 0\},
\end{align}
for $1\le i\le N$. We obtain the useful representation
\begin{align}\label{lambdarep}
\Lambda_t^N=\frac{\alpha}{N}\sum_{i=1}^N \mathbf{1}_{[0,t]}(\sigma^{N,i})-\frac{\beta}{N}\sum_{i=1}^N \mathbf{1}_{[0,t]}(\tau^{N,i}).
\end{align}

By the respective exchangeability of the particles $X^{N,1},\dots,X^{N,N}$ and $Y^{N,1},\dots, Y^{N,N}$, we see that this approximates $t\mapsto \alpha \bP(\sigma^{N,1}\le t)-\beta \bP(\tau^{N,1} \le t)$ with empirical averages.~We can now illustrate whence the physical jump condition~\eqref{eq:physcond} arises.~Suppose that there is a positive jump $\Lambda_t^N-\Lambda_{t-}^N>0$ at time $t$ for a particle system governed by the equations \eqref{eq:partsys}, \eqref{eq:sigmadef} and \eqref{lambdarep}. Consider the empirical measure $\rho^N_{t-}=\frac{1}{N}\sum_{i=1}^N \delta_{X_{t-}^{N,i}}\ind_{\{\sigma^{N,i}\ge t\}}$ of the $X^{N,i}$--particles that have not crossed zero by time $t-$. We see that any $x>0$ which solves $\alpha \rho^N_{t-}([0,x])=x$ gives a different solution to \eqref{eq:partsys}, \eqref{eq:sigmadef}, \eqref{lambdarep} with jump size $\Lambda_t^N-\Lambda_{t-}^N=x$. The physical condition enforces that one chooses the solution with the minimal jump size, and it is easy to check that this holds for the solutions we have constructed here (cf., for instance, \cite[Lemma~3.5]{cuchiero_propagation_2020} on the equivalence of physical and minimal solutions).

\begin{rmk}\label{rmk:1phasecompare}

When $\beta=0$, the subsystem $((X^{N,i})_{1\le i\le N},\Lambda^{N})$ in the construction above agrees with a physical solution to the one-phase problem
\begin{equation}\label{1phasecomp}
\begin{split}
& X_{t}^{i}=X_{0-}^{i}+B_t^{i}-\Lambda^{N}_t,\quad 1\le i\le N,\\
& \sigma^{i}=\inf\{t\ge 0: X_t^{i}\le 0\},\quad 1\le i\le N,\\
& \Lambda^{N}_t=\frac{\alpha}{N}\sum_{i=1}^N \mathbf{1}_{[0,t]}(\sigma^{i}),
\end{split}
\end{equation}
as studied in \cite{cuchiero_propagation_2020} (see also the similar systems in \cite{delarue_particle_2015}, \cite{nadtochiy_particle_2019}). In the proof of Theorem \ref{thm:exist}, we compare our two-phase particle system with $\beta>0$ to such a physical one-phase system.
\end{rmk}

\begin{rmk}\label{rmk:comparison}
It is also worth noting that the above construction of physical solutions was carried out pathwise; therefore, the driving Brownian motions may be replaced by arbitrary c\`adl\`ag paths $(G^{i})_{1\le i\le N}$. Consider the system
\begin{equation}
\begin{split}
&\widehat{X}_{t}^{i}=X_{0-}^{i}+G_t^{i}-\widehat{\Lambda}^{N}_t,\quad 1\le i\le N,\\
&\widehat{\sigma}^{i}=\inf\{t\ge 0:\,\widehat{X}_t^{i}\le 0\},\quad 1\le i\le N,\\
&\widehat{\Lambda}^{N}_t=\frac{\alpha}{N}\sum_{i=1}^N \mathbf{1}_{[0,t]}(\wh{\sigma}^{i}).
\end{split}
\end{equation}
There is a comparison principle between the physical solutions for this system and for~\eqref{1phasecomp}. If $G^{i}_t\le B_t^{i}$, $t\ge 0$, $1\le i\le N$, then $\widehat{\sigma}^{i}\le\sigma^{i}$, $1\le i\le N$ and $\widehat{\Lambda}^{N}_t\ge \Lambda^{N}_t$, $t\ge 0$.
\end{rmk}

\subsection{Application: A Systemic Risk Model With Two Competing Regions}\label{subsec:app}

We let $X^1,\dots,X^N$ denote the asset values of $N$ firms belonging to a common region of the world; we call these $\mathscr{X}$-type firms. Similarly, we take $N$ $\mathscr{Y}$-type firms with asset values $\wt{Y}^1,\dots,\wt{Y}^N$ from a competing region, and consider the dynamics
\begin{equation}\label{eq:tworegion}
\left\{
\begin{aligned}
    X_t^{i}&=X_{0-}^i+B_t^i-\Lambda_t^N,\quad t\ge 0,\\
    \wt{Y}_t^{i}&=-Y_{0-}^i-W_t^i+\Lambda_t^N,\quad t\ge 0,\\
    \sigma^i &= \inf \{t\ge 0: X_t^i\le 0\},\\
    \tau^i &= \inf \{t\ge 0: \wt{Y}_t^i\le 0\},\\
    \Lambda_t^N&=\frac{\alpha}{N}\sum_{i=1}^N \ind_{\{\sigma^i\le t\}}-\frac{\beta}{N}\sum_{i=1}^N \ind_{\{\tau^i\le t\}} ,\quad t\ge 0,
\end{aligned}
\right.
\end{equation}
which are obtained by applying the change of coordinates $\wt{Y}_t=-Y_t$ to \eqref{eq:partsys}, \eqref{eq:sigmadef}, and \eqref{lambdarep}. $\Lambda^N$ is taken to be the physical solution of this system, as developed above.
The assets of each firm follow a one-dimensional Brownian motion, as in Bachelier's model, except at the times $(\sigma^i)_{1\le i\le N}$ and $(\tau^i)_{1\le i\le N}$ when there are jumps in the interaction term $\Lambda^N$.
At time $\sigma^1$, for example, firm $X^1$ defaults by hitting zero, and in a ``flight-to-quality'' investors reposition wealth from region $\scrX$ to region $\scrY$, causing a loss of size $\alpha/N$ for all other $\scrX$-type firms and a gain of size $\alpha/N$ for all $\scrY$-type firms. If there are $\scrX$-type firms close to default at such a time, they may be pulled under, triggering a cascade of defaults. In total, the magnitude of the jump in the interaction term $\Lambda$ at time $\sigma^1$ is equal to $\alpha$ times the proportion of $\scrX$-type firms which default at that time.
The default of a $\scrY$-type firm at some time $\tau^i$ has a similar negative effect on $\scrY$-type firms and a beneficial effect on $\scrX$-type firms, with each defaulted firm contributing $\beta/N$ to the magnitude of the interaction term. 

\medskip

\begin{figure}
    \centering
    \hspace{0.7cm}
    \includegraphics[height=9.5cm]{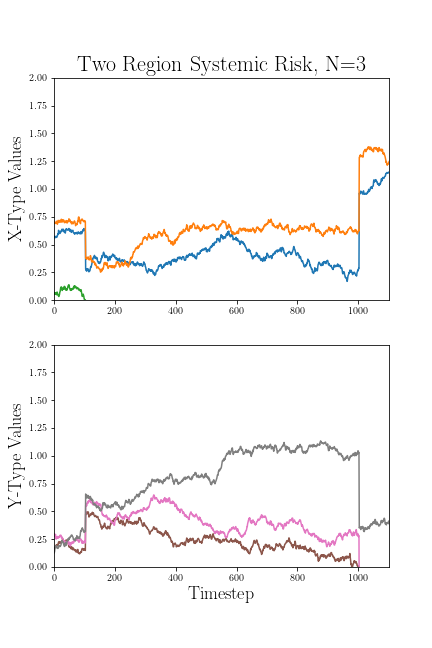}
    \includegraphics[height=9.5cm]{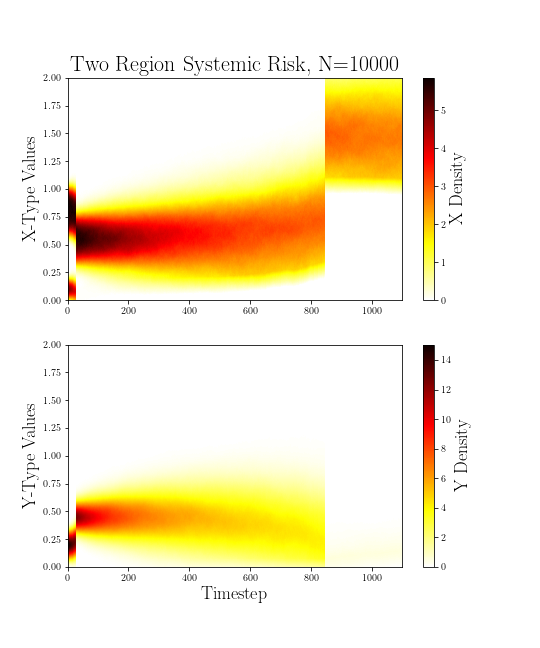}
    \caption{Realizations of \eqref{eq:tworegion} with initial densities given by \eqref{eq:exampledenisties}.}
    \label{fig:plots}
\end{figure}

To illustrate this model, suppose that the $\mathscr{X}$-type firms are banks in the United States and the $\mathscr{Y}$-type firms are banks in Canada. Let $\alpha=\beta=1$ and let $X_{0-}$ and $Y_{0-}$ have densities 
\begin{align}\label{eq:exampledenisties}
f=\frac{5}{2}\cdot\ind_{[0.05,0.15]}+\frac{15}{8}\cdot\ind_{[0.60,1.00]}\quad\text{and}\quad g=5\cdot \ind_{[0.10,0.30]},
\end{align}
respectively. In this scenario, a quarter of the US banks are initially close to default with asset values in the interval $[0.05,0.15]$, and the remainder of the US banks start in $[0.60,1.00]$. All of the Canadian banks are initially in the moderate interval $[0.10,0.30]$.
In Figure \ref{fig:plots}, we plot realizations of this system when $N=3$ and when $N=10000$. For $N=3$, each sample path is pictured up until the respective bank's time of default; we refer to these as absorbed paths. For $N=10000$, we plot the densities of the absorbed paths.
For these realizations, we see that the subset of US banks which start in $[0.05,0.15]$ default early on, giving a boost to the Canadian banks as investors move their wealth to that region. However, the boost is not enough to prevent a systemic event at a later time when a majority of the Canadian banks default simultaneously. With Theorems \ref{tight} and \ref{thm:exist}, we will see that the large system limit as $N\to\infty$ satisfies the McKean--Vlasov system \eqref{limitprob}. We will return to this example in Section \ref{sec:conc} to discuss the implications for systemic risk and to indicate possible extensions.

\begin{rmk}\label{rmk:alphabeta}
The constants $\alpha>0$ and $\beta>0$ in the above model measure the strength of self- and cross-excitation for type $\scrX$ and $\scrY$ firms, respectively.  The coupling between the two regions results in the possible non-monotonicity of the interaction term $\Lambda^N$, which places the current model beyond the mathematical framework of the one-phase financial contagion model in \cite{nadtochiy_particle_2019}.
As mentioned in the introduction, the case $\alpha \ge 0 \ge \beta$ (corresponding to a supercooled liquid meeting a non-superheated solid in the language of phase change problems) results in a monotone boundary, which simplifies the proof of tightness for the laws of the finite particle systems (for instance, compare \cite[proof of Lemma 5.4]{delarue_particle_2015} to the proof of the result in the next section).
\end{rmk}

\section{Tightness in M1 Topology}\label{sec:tight}

To obtain a large system limit we follow the usual compactness approach: first, we prove the tightness of a family of random processes (in a suitable topology), then we extract the limit of a convergent subsequence and show it solves \eqref{limitprob}. Our objects of study are the random empirical measures $(\mu^N,\nu^N):=\big(\frac{1}{N}\sum_{i=1}^N \delta_{X^{N,i}},\frac{1}{N}\sum_{i=1}^N \delta_{Y^{N,i}}\big)$. We consider $\mu^N$ and $\nu^N$ as probability measures on path space, endowed with the M1 topology which we describe next.

\subsection{M1 Topology}\label{subsec:M1}

For any given $S,T\in\R$ with $S<T$, let $\widetilde{D}([S,T],\R)$ denote the space of c\`adl\`ag functions from $[S,T]$ to $\R$ that are left-continuous at $T$. In the sequel, we will also let $C([S,T],\R)$ denote the space of continuous functions from $[S,T]$ to $\R$.
We now give some basic facts about the M1 topology on $\widetilde{D}([S,T],\R)$; for a comprehensive introduction we refer the reader to \cite[Chapter 12]{whitt_stochastic-process_2002}. For $h\in \widetilde{D}([S,T],\R)$, we define the \emph{completed graph} of $h$:
\begin{align}
\calG_h:=\big\{(x,t)\in\R\times [S,T]: x\in [h(t-)\wedge h(t),h(t-)\vee h(t)]\big\}.
\end{align}
$\calG_h$ is equipped with the total ordering $\le$, where $(x_1,t_1)\le(x_2,t_2)$ if (i) $t_1< t_2$ or (ii) $t_1=t_2$ and $|h(t_1-)-x_1|\le|h(t_1-)-x_2|$ (this is the natural order when we trace the completed graph $\calG_h$ from left to right). A \emph{parametric representation} of $h$ is a continuous function $(u,r)$ mapping $[0,1]$ onto $\calG_h$ which is non-decreasing with respect to the order $\le$. We let $\calR_h$ denote the set of all parametric representations for a given $h$ and define the M1 distance between $h_1,h_2\in \widetilde{D}([S,T],\R)$ by
\begin{align}
d_{\text{M1}}(h_1,h_2)=\inf_{(u^{h_1},r^{h_1})\in\calR_{h_1}, (u^{h_2},r^{h_2})\in\calR_{h_2}}\lb(\|u^{h_1}-u^{h_2}\|_\infty\vee\|r^{h_1}-r^{h_2}\|_\infty\rb).
\end{align}
This defines a metric which induces the M1 topology on $\widetilde{D}([S,T],\R)$. We note here that $\widetilde{D}([S,T],\R)$ with M1 is Polish (see \cite[Section 12.8]{whitt_stochastic-process_2002}).

\medskip

We now introduce some oscillation functions that help to characterize convergence in M1. For $h\in \widetilde{D}([S,T],\R)$, $\delta>0$, and $t\in[S,T]$, we define the $w$--\emph{oscillation} by
\begin{align}
w_{[S,T]}(h,t,\delta)=\sup_{S\vee(t-\delta)\leq t_1 < t_2 < t_3\le T\wedge (t+\delta)}\inf_{\theta\in[0,1]} \big|\theta h(t_1)+(1-\theta)h(t_3)-h(t_2)\big|.
\end{align}
Note that
\begin{equation}
\begin{split}
&\inf_{\theta\in[0,1]} \big|\theta h(t_1)+(1-\theta)h(t_3)-h(t_2)\big| \\
&=\begin{cases}
0,&\text{if}\;\;h(t_2)\in[h(t_1),h(t_3)], \\
|h(t_2)-h(t_1)|\wedge |h(t_2)-h(t_3)|,&\text{otherwise}
\end{cases}
\end{split}
\end{equation}
is the distance between $h(t_2)$ and the line segment $[h(t_1),h(t_3)]$. We write
\begin{align}
v_{[S,T]}(h,t,\delta):=\sup_{S\vee(t-\delta)\leq t_1 \le t_2\le T\wedge (t+\delta)}|h(t_2)-h(t_1)|
\end{align}
for the modulus of continuity around a point $t\in[S,T]$.~Finally, we put $w_{[S,T]}$ and $v_{[S,T]}$ together to give the M1--\emph{oscillation}
\begin{align}
u_{[S,T]}(h,\delta):=v_{[S,T]}(h,S,\delta)\vee v_{[S,T]}(h,T,\delta)  \vee \sup_{t\in[S,T]}w_{[S,T]}(h,t,\delta).
\end{align}

We recall the following tightness criterion.

\begin{lem}[Theorem 12.12.3 in \cite{whitt_stochastic-process_2002}]\label{criterion}
A sequence of probability measures $(\bP_N)_{N\ge 1}$ on $\widetilde{D}([S,T],\R)$, endowed with the M1 topology, is tight if and only if both
\begin{enumerate}
\item for each $\veps>0$, there exists some $C<\infty$ such that
\begin{align*}
\bP_N\big(\big\{h\in\widetilde{D}([S,T],\R):\,\|h\|_\infty>C\big\}\big)<\veps,\quad N\ge 1,
\end{align*}
and
\item for each $\veps>0$ and $r>0$, there exists some $\delta>0$ such that
\begin{align*}
\bP_N\big(\big\{h\in\widetilde{D}([S,T],\R):\, u_{[S,T]}(h,\delta) \ge r\big\}\big)<\veps,\quad N\ge 1.
\end{align*}
\end{enumerate}
\end{lem}

\subsection{Tightness}\label{subsec:tight}
In Theorem \ref{tight} below, we show the tightness of the sequence $(\mu^N,\nu^N)_{N\ge 1}$ of random empirical measures. Noting the $v_{[S,T]}(h,S,\delta)$ and $v_{[S,T]}(h,T,\delta)$ terms in $u_{[S,T]}(h,\delta)$, we must ensure the continuity of the particles' sample paths at the endpoints of our time domain. To do so, we embed our processes into $\widetilde{D}([-1,T+1],\R)$:
\begin{equation}\label{extend}
\begin{split}
&B_t^{i}:=0\text{ for }t\in[-1,0)\text{ and }W_t^{i}:=0\text{ for }t\in[-1,0)\text{, for all }1\le i\le N;\\
&X_t^{N,i}:=X_{0-}^{N,i}\text{ for }t\in[-1,0)\text{ and }Y_t^{N,i}:=Y_{0-}^{N,i}\text{ for }t\in[-1,0)\text{, for all }1\le i\le N;\\
&\text{and }\Lambda^N_t:=0 \text{ for } t\in [-1,0)\text{ and }\Lambda^N_t:=\Lambda^N_T \text{ for } t\in (T,T+1].
\end{split}
\end{equation}
In other words, we make the particles' sample paths constant to the left of $0$, and the interaction term constant to the left of 0 and to the right of $T$. As observed in \cite[Remark 4.5]{delarue_particle_2015}, the extension to the right is without loss of generality since we can construct solutions on a sequence of intervals $[-1,T_n+1]$ with $T_n\to\infty$ and use Cantor's diagonal argument to get a solution on $[-1,\infty)$. We use the extension to the left since under (the fairly weak) Assumption \ref{ass1} we must allow for a discontinuity at~0 in the limiting interaction term $\Lambda$.
The next lemma bounds the $w$--oscillation of~$\Lambda^N$.
\begin{lem}\label{lem:w_T}
For any $N\ge 1$, $t\in[0,T]$, $\delta>0$, and $r>0$,
\begin{align}
\bP\big(w_{[0,T]}(\Lambda^N,t,\delta)\ge r\big)\le \frac{4(\alpha+\beta)}{r}\,\Phi\lb(-\frac{r}{2\sqrt{2\delta}}\rb).
\end{align}
\end{lem}

\noindent\textbf{Proof.} From the definition of $w_{[0,T]}$, we see that $w_{[0,T]}(\Lambda^N,t,\delta)\ge r$ only if for some $t_2\in[0\vee(t-\delta),T\wedge (t+\delta)]$, there exist $t_1\in[0\vee(t-\delta),t_2)$ and $t_3\in(t_2,T\wedge (t+\delta)]$ such that $|\Lambda^N_{t_2}-\Lambda^N_{t_1}|\ge r$, $|\Lambda^N_{t_3}-\Lambda^N_{t_2}|\ge r$ and $\sgn(\Lambda^N_{t_2}-\Lambda^N_{t_1})\neq \sgn(\Lambda^N_{t_3}-\Lambda^N_{t_2})$. 

\medskip

We first consider the case where $\Lambda^N_{t_2}-\Lambda^N_{t_1}\le -r$ and $\Lambda^N_{t_3}-\Lambda^N_{t_2}\ge r$, which we denote by $E_{lr}$, the ``left-right'' event. The other event (``right-left'') is denoted by $E_{rl}$. Let
\begin{eqnarray}
&& \eta_1 := \inf\big\{s>0\vee(t-\delta):\,\Lambda^N_{s}-\Lambda^N_{s'}\le-r\text{ for some }s'\in[0\vee(t-\delta),s)\big\}, \\
&&\eta_2 := \inf\bigg\{s>\eta_1:\,\frac{\alpha}{N}\sum_{i=1}^N \mathbf{1}_{(\eta_1,s]}(\sigma^{N,i})\ge \frac{r}{2}\bigg\}.
\end{eqnarray}
Additionally, we let 
\begin{align}
	\eta_0 := \inf\big\{s'\in[0\vee(t-\delta),\eta_1):\,\Lambda^N_{\eta_1}-\Lambda^N_{s'}\le-r\big\},
\end{align}
and note that $\Lambda^N_{\eta_1}-\Lambda^N_{\eta_0}\le-r$ due to the right-continuity of $\Lambda^N$.
It is clear that $E_{lr}$ implies the event $\{0\vee(t-\delta)\le \eta_0<\eta_1<\eta_2\le T\wedge (t+\delta)\}$.
By Markov's inequality and exchangeability, we obtain the upper bound
\begin{equation}
	\begin{split}
		\bP(E_{lr})
		&\le \bP\bigg(\frac{\alpha}{N}\sum_{i=1}^N \mathbf{1}_{(\eta_1,T\wedge(t+\delta)\wedge\eta_2]}(\sigma^{N,i})\ge \frac{r}{2},\;\eta_0<\eta_1< T\wedge (t+\delta)\bigg) \\
		&\le\frac{2}{r}\,\E\bigg[\frac{\alpha}{N}\sum_{i=1}^N \mathbf{1}_{(\eta_1,T\wedge(t+\delta)\wedge\eta_2]}(\sigma^{N,i})\,\mathbf{1}_{\{\eta_0<\eta_1< T\wedge (t+\delta)\}}\bigg]\\
		&=\frac{2 \alpha}{r}\,\bP\lb(\sigma^{N,1}\in(\eta_1,T\wedge (t+\delta)\wedge\eta_2],\;\eta_0<\eta_1< T\wedge (t+\delta)\rb).
	\end{split}
\end{equation}

\smallskip

Next, ``ignoring'' the presence of the $X^{N,1}$--particle, we define
\begin{eqnarray}
	&& \wt{\Lambda}^N_t=\frac{\alpha}{N}\sum_{i=2}^N \mathbf{1}_{[0,t]}(\widetilde{\sigma}^{N,i})-\frac{\beta}{N}\sum_{i=1}^N \mathbf{1}_{[0,t]}(\widetilde{\tau}^{N,i}), \\
	&& \wt\eta_1 = \inf\big\{s>0\vee(t-\delta):\,\wt\Lambda^N_{s}-\wt\Lambda^N_{s'}\le-r\text{ for some }s'\in[0\vee(t-\delta),s)\big\},  \\
	&&\wt\eta_0 = \inf\big\{s'\in[0\vee(t-\delta),\wt\eta_1):\,\wt\Lambda^N_{\wt\eta_1}-\wt\Lambda^N_{s'}\le-r\big\}.
\end{eqnarray}
On the event $\{\sigma^{N,1}\in(\eta_1,T\wedge (t+\delta)\wedge\eta_2],\,\eta_0<\eta_1< T\wedge (t+\delta)\}$, we note that $\wt\eta_1=\eta_1$ and $\wt\eta_0=\eta_0$ since $\wt{\Lambda}^N_t=\Lambda^N_t$ before time $\sigma^{N,1}$.  Hence,
\begin{equation}
	\begin{split}
		&\;\bP\lb(\sigma^{N,1}\in(\eta_1,T\wedge (t+\delta)\wedge\eta_2],\;\eta_0<\eta_1< T\wedge (t+\delta)\rb)\\
		&\le\bP\lb(\sigma^{N,1}\in(\wt\eta_1,T\wedge (t+\delta)\wedge\eta_2]\,\big|\,\wt\eta_0<\wt\eta_1< T\wedge (t+\delta)\rb)\bP\lb(\wt\eta_0<\wt\eta_1< T\wedge (t+\delta)\rb)\\
		&\le\bP\lb(\sigma^{N,1}\in(\wt\eta_1,T\wedge (t+\delta)\wedge\eta_2]\,\big|\,\wt\eta_0<\wt\eta_1< T\wedge (t+\delta)\rb).
	\end{split}
\end{equation}
For $\wt{\eta_1}\le s\le \eta_2$, we have the comparisons
\begin{align}
	X_{\wt\eta_1}^{N,1}&=X_{\wt\eta_0}^{N,1} + (B_{\wt\eta_1}^{1}-B_{\wt\eta_0}^{1})-(\wt\Lambda_{\wt\eta_1}^N-\wt\Lambda_{\wt\eta_0}^N)\ge X_{\wt\eta_0}^{N,1} + (B_{\wt\eta_1}^{1}-B_{\wt\eta_0}^{1})+r
\end{align}
and
\begin{equation}
	\begin{split}
		X_{\wt{s}}^{N,1}&=X_{\wt\eta_1}^{N,1} + (B_{s}^{1}-B_{\wt\eta_1}^{1})-(\Lambda_{s}^N-\Lambda_{\wt\eta_1}^N)\\
		&\ge \big(X_{\wt\eta_0}^{N,1} + (B_{\wt\eta_1}^{1}-B_{\wt\eta_0}^{1})+r\big) + (B_{s}^{1}-B_{\wt\eta_1}^{1})-(\Lambda_{s}^N-\Lambda_{\wt\eta_1}^N)\\
		&\ge  X_{\wt\eta_0}^{N,1} + (B_{s}^{1}-B_{\wt\eta_0}^{1})+\frac{r}{2}\label{Xcomp},
	\end{split}
\end{equation}
where we have used the definition of $\eta_2$ to bound $\Lambda_{\wt{s}}^N-\Lambda_{\wt\eta_1}^N$ by $\frac{r}{2}$. We also must have $X_{\wt\eta_0}^{N,1}>0$ on the event $\{\sigma^{N,1}\in (\wt\eta_1,T\wedge (t+\delta)\wedge\eta_2]\}$. Our comparison \eqref{Xcomp} gives
\begin{equation}
	\begin{split}
		&\;\bP\lb(\sigma^{N,1}\in (\wt\eta_1,T\wedge (t+\delta)\wedge\eta_2]\,\big|\,\wt\eta_0<\wt\eta_1< T\wedge (t+\delta)\rb)\\
		&\le  \bP\Big(\inf_{\wt\eta_1\le s\le T\wedge(t+\delta)\wedge\eta_2}\Big( (B_{s}^{1}-B_{\wt\eta_0}^{1})+\frac{r}{2}\Big)\le 0\,\Big|\,\wt\eta_0<\wt\eta_1< T\wedge (t+\delta)\Big) \\
		&\le  \bP\Big(\inf_{\widetilde{\eta}_0\le s\le T\wedge(t+\delta)}(B_{s}^{1}-B_{\wt\eta_0}^{1})\le -\frac{r}{2}\,\Big|\,\wt\eta_0<\wt\eta_1< T\wedge (t+\delta)\Big) 
		\le\bP\Big(\inf_{0\le s\le 2\delta} B_{s}^{1}\le -\frac{r}{2}\Big),
	\end{split}
\end{equation}
where the final inequality follows from the independence of $B^1$ from $\wt\Lambda^N$, $\wt\eta_0$ and $\wt\eta_1$. Putting things together and using the reflection principle for Brownian motion yields
\begin{equation}
		\bP\lb(E_{lr}\rb)\le \frac{2\alpha}{r}\,\bP\lb(|B_{2\delta}^{1}|\ge\frac{r}{2}\rb)\\
		=\frac{4\alpha}{r}\,\Phi\lb(-\frac{r}{2\sqrt{2\delta}}\rb).
\end{equation}
With a similar estimate for $\bP\lb(E_{rl}\rb)$, we obtain
\begin{equation}
\quad\bP\lb(w_{[0,T]}(\Lambda^N,t,\delta)\ge r\rb)
\le \bP\lb(E_{rl}\rb)+\bP\lb(E_{lr}\rb)\\
\le \frac{4(\alpha+\beta)}{r}\,\Phi\lb(-\frac{r}{2\sqrt{2\delta}}\rb).\qquad\qed
\end{equation}

\smallskip

We now give our tightness result.~For $N\ge 1$, we consider the empirical measures $\mu^N$ and $\nu^N$ as $\calP(\widetilde{D}([-1,T+1],\R))$--valued random variables, where the particles' sample paths have been extended as described at the beginning of this section.~Here, $\calP(\widetilde{D}([-1,T+1],\R))$ denotes the space of probability measures on $\widetilde{D}([-1,T+1],\R)$, endowed with the topology of weak convergence induced by the M1 topology.~To clarify terminology, a sequence of random variables is said to be \emph{tight} if the sequence of their probability laws is tight. 
\begin{thm}\label{tight}
Let Assumption \ref{ass1} hold.~Then, $(\mu^N,\nu^N,\Lambda^N)_{N\ge 1}$ is tight on the space  $\calP(\widetilde{D}([-1,T+1],\R))\times \calP(\widetilde{D}([-1,T+1],\R))\times \widetilde{D}([-1,T+1],\R)$.
\end{thm}

\noindent\textbf{Proof.} We only show that $(\mu^N)_{N\ge 1}$ is tight on $\calP(\widetilde{D}([-1,T+1],\R))$, since the other components can be analyzed very similarly. 
~Using a standard result from the theory of propagation of chaos (\cite[Proposition 2.2]{burkholder_topics_1991}, which applies since~$\widetilde{D}([-1,T\!+\!1],\R)$ with the M1 topology is Polish), it suffices to prove that $(X^{N,1})_{N\ge 1}$ is tight on $\widetilde{D}([-1,T+1],\R)$.~We recall that $X_{0-}^{N,1}\stackrel{d}{=}\mu_{0-}$, with $\mu_{0-}$ having a finite first moment by Assumption \ref{ass1}, and that $X_t^{N,1}=X_{0-}^{N,1}$ for $t\in[-1,0)$ by construction.

\medskip

We aim to apply the tightness criterion of Lemma \ref{criterion}.~From the representation \eqref{lambdarep}, we see that $|\Lambda_t^N|\le \alpha \vee \beta$ for all $t\in [-1,T+1]$ and $N\ge 1$.~Hence,
\begin{equation}
\begin{split}
\bP\big(\sup_{t\in[-1,T+1]}|X_t^{N,1}|> C\big)&=\bP\big(\sup_{t\in[0,T+1]}|X_{0-}^{N,1}+B_t^{1}-\Lambda_t^N|> C\big)\\
&\le \bP\big(X_{0-}^{N,1}+\sup_{t\in[0,T+1]}|B_t^{1}|> C-\alpha\vee\beta\big),
\end{split}
\end{equation}
which can be made arbitrarily small, independently of $N$, by taking $C$ large enough. This confirms the boundedness part of the tightness criterion in Lemma \ref{criterion}.

\medskip

For the M1--oscillation part of the tightness criterion, we estimate
\begin{equation}
\begin{split}
&\;\bP\big(u_{[-1,T+1]}(X^{N,1},\delta) \ge r\big)\\
&\le \bP\big(\sup_{t\in[-1,T+1]}w_{[-1,T+1]}(X^{N,1},t,\delta)\ge r\big) +\bP\big(v_{[-1,T+1]}(X^{N,1},-1,\delta)\ge r\big)\\
&\quad +\bP\big(v_{[-1,T+1]}(X^{N,1},T+1,\delta)\ge r\big).
\end{split}
\end{equation}
Since $X^{N,1}$ is constant on $[-1,0)$, we have $\bP(v_{[-1,T+1]}(X^{N,1},-1,\delta)\!\ge\! r)\!=\!0$, $\delta
\!\in\!(0,1)$. On $[T,T+1]$, the process $X^{N,1}$ has the modulus of continuity of the Brownian motion $B^1$, hence $\bP\big(v_{[-1,T+1]}(X^{N,1},T+1,\delta)\ge r\big)$ can be made arbitrarily small, independently of $N$, by taking $\delta>0$ sufficiently small.

\medskip

Since $\mu_{0-}$ and $\nu_{0-}$ do not have atoms at zero, we have $\sigma^{N,i}>0$ and $\tau^{N,i}>0$ for all $1\le i\le N$ and $N\ge 1$, and thus $\Lambda_0^{N}=\Lambda_{0-}^{N}$.~Following \cite[proof of Corollary~12.7.1]{whitt_stochastic-process_2002} and since $\Lambda^N$ is constant on $[-1,0]$ and $[T,T+1]$, we have
\begin{equation}
\begin{split}
&\;\bP\big(\sup_{t\in[-1,T+1]}w_{[-1,T+1]}(X^{N,1},t,\delta)\ge r\big) \\
&\le \bP\big(\sup_{t\in[0,T+1]}\big(w_{[0,T+1]}(\Lambda^{N},t,\delta)+2v_{[0,T+1]}(B^{1},t,\delta)\big)\ge r\big)\\
&\le \bP\Big(\sup_{t\in[0,T]}w_{[0,T]}(\Lambda^{N},t,\delta) \ge \frac{r}{2}\Big)+\bP\Big(\sup_{t\in[0,T+1]} v_{[0,T+1]}(B^{1},t,\delta)\ge \frac{r}{4}\Big).
\end{split}
\end{equation}
Dividing $[0,T]$ into $\lceil T/2\delta \rceil$ overlapping intervals of length $4\delta$, taking a union bound, and applying Lemma \ref{lem:w_T} on each interval gives
\begin{equation}\label{wTestimate}
\begin{split}
\bP\Big(\sup_{t\in[0,T]}w_{[0,T]}(\Lambda^{N},t,\delta)\ge \frac{r}{2}\Big) &\le \sum_{k=1}^{\lceil T/2\delta \rceil}\bP\Big(w_{[0,T]}(\Lambda^{N},2k\delta,2\delta)\ge \frac{r}{2}\Big) \\
&\le \lceil T/2\delta \rceil\,\frac{8(\alpha+\beta)}{r}\,\Phi\lb(-\frac{r}{8\sqrt{\delta}}\rb).
\end{split}
\end{equation}
Using the estimate $\Phi(-x)\le \frac{\varphi(x)}{x}$ (that can be obtained through integration by parts) and since $\lim_{\delta\downarrow0} \delta^{-1/2}\,e^{-r^2/(128\delta)}=0$, for any given $r,\veps>0$ we can find some~$\delta>0$ such that $\bP\lb(\sup_{t\in[0,T]}w_{[0,T]}(\Lambda^{N},t,\delta)\ge \frac{r}{2}\rb)<\frac{\veps}{3}$. By making use of the almost sure continuity of Brownian sample paths, and possibly after taking $\delta$ smaller, we can ensure that $\bP\lb(\sup_{t\in[0,T+1]} v_{[0,T+1]}(B^{1},t,\delta)\ge \frac{r}{4}\rb)<\frac{\veps}{3}$.
Putting the estimates together we see that the tightness criterion of Lemma \ref{criterion} is verified. \qed

\section{Large System Limit}\label{sec:limit}

We introduce the empirical measures
\begin{align}
\zeta^N:=\frac{1}{N}\sum_{i=1}^N\delta_{(X^{N,i},B^{i})}\quad\text{and}\quad\theta^N:=\frac{1}{N}\sum_{i=1}^N\delta_{(Y^{N,i},W^{i})}
\end{align}
on $\widetilde{D}([-1,T+1],\R)\times C([-1,T+1],\R)$ which have the respective first marginals $\mu^N$ and $\nu^N$.~Here, again, we extend our processes to $[-1,0)$ and $(T,T+1]$ via \eqref{extend}.~Let $(\Pi^N)_{N\geq1}$ denote the laws of $(\zeta^N,\theta^N,\Lambda^N)_{N\ge 1}$,
which are tight thanks to Theorem \ref{tight}. Consider a limit point $\Pi$ of some convergent subsequence of $(\Pi^N)_{N\geq1}$, where we choose not to relabel the subsequence.~For $(\zeta,\theta,\Lambda)$ sampled from $\Pi$, let $\mu$ and $\nu$ denote the first marginals of $\zeta$ and $\theta$, respectively. Note that $\mu^N\to\mu$ and $\nu^N\to \nu$ weakly in law by the Continuous Mapping Theorem.

\medskip

Furthermore, for $t\in[0,T]$, we consider the functions $m_t:\,\widetilde{D}([-1,T+1],\R)\to \R$ and $n_t:\,\widetilde{D}([-1,T+1],\R)\to \R$ defined by 
\begin{equation}
m_t(h_1)=\ind_{\{\inf_{s\in[0,t]} h_1(s) \le 0\}}\quad\text{and}\quad
n_t(h_2)=\ind_{\{\sup_{s\in[0,t]} h_2(s) \ge 0\}}.
\end{equation}
Equipped with this notation and under Assumption \ref{ass1} we have:

\begin{thm}\label{thm:exist}
Let $(x,b,y,w)$ denote the canonical process on $\widetilde{D}([-1,T+1],\R)\times C([-1,T+1],\R)\times \widetilde{D}([-1,T+1],\R)\times C([-1,T+1],\R)$.~Under $\Pi$--almost every $(\zeta,\theta,\Lambda)$, the quintuple $(x_t,b_t,y_t,w_t,\Lambda_t)_{t\in[-1,T+1]}$ with $(x_t,b_t)$ drawn from $\zeta$ and $(y_t,w_t)$ drawn from $\theta$ generates a physical solution of \eqref{limitprob} on $[0,T)$.~That is,
\begin{enumerate}[label=(\roman*)]
\item \label{thm:bm}$x\!=\!b\!-\!\Lambda$ and $y\!=\!w\!-\!\Lambda$, where $b$ and $w$ are standard Brownian motions on $[0,T)$.
\item \label{thm:init}$\Lambda_{0-}\!=\!0$, $x_{0-}$ is distributed according to $\mu_{0-}$, and $y_{0-}$ according to~$\nu_{0-}$.
\item $\Lambda_t=\alpha\E^\mu[m_t]-\beta\E^\nu [n_t]$ for $t\in[0,T)$. \label{thm:lambda}
\item \label{thm:phys} $\Lambda$ is physical on $[0,T)$ in the sense of Definition \ref{def:phys}.
\end{enumerate}
\end{thm}

The proof of \ref{thm:lambda} is based on \cite[proof of Theorem 4.4]{delarue_particle_2015} (see also \cite[proof of Theorem 2.4]{nadtochiy_particle_2019}).~We are able to simplify some of the steps by keeping track of the Brownian motions in $\zeta^N$ and $\theta^N$. Playing the role of the crossing property which is used in prior works (see, e.g., \cite[Lemma 5.6]{delarue_particle_2015}), we use the following stronger property of Brownian sample paths.

\begin{lem}[Crossing Property]\label{lem:cross}
Let $(B_t)_{t\ge 0}$ be a standard Brownian motion, let $(h(t))_{t\ge 0}$ be a c\`adl\`ag function, and define $X=B-h$. Almost surely, $X$ does not have any zeroes in $(0,\infty)$ which are local extrema.
\end{lem}

\noindent\textbf{Proof.}
The proof of \cite[Proposition 3.4(i)]{antunovic_isolated_2011} applies. \qed

\begin{rmk}\label{rmk:cross}
The expression for $\Lambda$ in Theorem \ref{thm:exist}\ref{thm:lambda} deserves some explanation. We recall the representation \eqref{lambdarep} for the finite particle system:
\begin{align*}
\Lambda_t^N=\frac{\alpha}{N}\sum_{i=1}^N\ind_{[0,t]}(\sigma^{N,i})-\frac{\beta}{N}\sum_{i=1}^N\ind_{[0,t]}(\tau^{N,i}).
\end{align*}
For each $1\le i\le N$ and $N\ge1$, the trajectory of the particle $X^{N,i}$ satisfies the crossing property almost surely, and hence the functions $t\mapsto\mathbf{1}_{[0,t]}(\sigma^{N,i})$ and $t\mapsto m_t(X^{N,i})$ are identical on $[0,T]$.~This holds similarly for the $Y^{N,i}$--particles, so
\begin{align}
\Lambda_t^N=\frac{\alpha}{N}\sum_{i=1}^Nm_t(X^{N,i})-\frac{\beta}{N}\sum_{i=1}^Nn_t(Y^{N,i})
=\alpha \E^{\mu^N}[m_t]-\beta \E^{\nu^N}[n_t],\quad t\in[0,T]
\end{align}
almost surely. For the limiting system, \eqref{limitprob} stipulates
\begin{align*}
\Lambda_t&=\alpha \bP(\sigma\le t)-\beta \bP(\tau\le t),
\end{align*}
where $\sigma$ and $\tau$ are the crossing times of zero for the representative particles $X$ and~$Y$, respectively.~By construction, $x$ and $y$ in Theorem \ref{thm:exist} have laws $\mu$ and $\nu$.~In the proof of the theorem we are going to see that, almost surely, $\mu$ and $\nu$ give mass only to trajectories satisfying the crossing property.~Hence,
\begin{align}
\Lambda_t=\alpha\E^\mu[m_t]-\beta\E^\nu [n_t]=\alpha \E^\mu[\ind_{[0,t]}(\sigma)]-\beta \E^\nu[\ind_{[0,t]}(\tau)],
\end{align}
as desired.
\end{rmk}

\noindent\textbf{Proof of Theorem \ref{thm:exist}.} We prove \ref{thm:bm} to \ref{thm:phys} in order. 

\medskip

\noindent\textbf{Proof of \ref{thm:bm}.}~Using the Skorokhod Representation Theorem (SRT) we can find $(\wt{\zeta}^N,\wt{\theta}^N,\wt{\Lambda}^N)_{N\ge 1}$ with laws $(\Pi^N)_{N\ge 1}$ converging almost surely to $(\wt{\zeta},\wt{\theta},\wt{\Lambda})$ with law~$\Pi$. By a second use of the SRT, we may take $(x^N,b^N)$ drawn from $\wt{\zeta}^N(\omega)$ and $(y^N,w^N)$ drawn from $\wt{\theta}^N(\omega)$ with $(x^N,b^N)\to (x,b)$ and $(y^N,w^N)\to (y,w)$ almost surely, where $(x,b)$ has law $\wt{\zeta}(\omega)$ and $(y,w)$ has law $\wt{\theta}(\omega)$. 
~From this setup, we must have that $b$ and~$w$ are standard Brownian motions on $[0,T+1]$, $\Pi$--almost surely. We note that $x^N=b^N-\widetilde{\Lambda}^N(\omega)$ holds for each $N\ge1$, with $b^N\to b$.~Since the $b$ summand in $b-\widetilde{\Lambda}(\omega)$ is a continuous function of time, by \cite[Theorem 12.7.3]{whitt_stochastic-process_2002} we have that $b^N-\widetilde{\Lambda}^N(\omega)\to b-\widetilde{\Lambda}(\omega)$. Using the same argument for $w-\widetilde{\Lambda}(\omega)$ gives \ref{thm:bm}.

\medskip

\noindent\textbf{Proof of \ref{thm:init}.} By construction, $\widetilde{\Lambda}^N_t=0$ for all $t\in[-1,0)$ and $N\ge1$.~Therefore, $\widetilde{\Lambda}_t=0$ for all $t\in[-1,0)$, and thus $\widetilde{\Lambda}_{0-}=0$.~Similarly, $x^N_t=x^N_{0-}$ for all $t\in[-1,0)$ and $N\ge1$, so that $x_t=x_{0-}$ for all $t\in[-1,0)$, $\widetilde{\zeta}(\omega)$--almost surely under $\Pi$.~Moreover, the law of $x_{0-}=x_{-1}$ under $\widetilde{\zeta}(\omega)=\lim_{N\to\infty} \widetilde{\zeta}^N(\omega)$ is given by $\mu_{0-}$, since by the Glivenko-Cantelli Theorem there are copies of $\frac{1}{N}\sum_{i=1}^N \delta_{X^{N,i}_{-1}}$, $N\ge1$ converging weakly to $\mu_{0-}$ almost surely.~Analogously, the law of $y_{0-}$ under $\widetilde{\theta}(\omega)$ is $\nu_{0-}$, $\Pi$--almost surely.

\medskip

\noindent\textbf{Proof of \ref{thm:lambda}.} In light of Remark \ref{rmk:cross}, we have $\widetilde{\Lambda}_t^N=\alpha \E^{\widetilde{\mu}^N}[m_t]-\beta \E^{\widetilde{\nu}^N}[n_t]$ for all $t\in[0,T]$, with $\widetilde{\Lambda}^N$ on the left-hand side of this equation converging to $\widetilde{\Lambda}$ in the M1 sense. Following \cite[proof of Theorem 4.4]{delarue_particle_2015} we can find a countable set $J\subset[0,T]$ such that for $\Pi$--almost every $\omega$, it holds $\E^{\widetilde{\mu}(\omega)}[m_t]=\E^{\widetilde{\mu}(\omega)}[m_{t-}]$ and $\E^{\widetilde{\nu}(\omega)}[n_t]=\E^{\widetilde{\nu}(\omega)}[n_{t-}]$ for $t\in[0,T]\backslash J$. For these $\omega$, we claim that 
\begin{align}\label{eq:mnconv}
\E^{\wt\mu^N(\omega)}[m_t]\to\E^{\wt\mu (\omega)}[m_t]\quad\text{and}\quad \E^{\wt\nu^N(\omega)}[n_t]\to\E^{\wt\nu(\omega)}[n_t],\quad t\in(0,T)\backslash J.
\end{align}
With minor modifications to \cite[Lemma 5.6 and Proposition 5.8]{delarue_particle_2015}, \eqref{eq:mnconv} holds if $(\wt{\zeta}(\omega),\wt{\theta}(\omega))$ gives full measure to the set of paths
\begin{align}
A=\big\{(x,b,y,w):\, x \text{ and } y \text{ have no zeroes on $(0,T)$ which are local extrema}\big\},
\end{align}
$\Pi$--almost surely.~Lemma \ref{lem:cross} for $x\!=\!b\!-\!\widetilde{\Lambda}(\omega)$ and $y\!=\!w\!-\!\widetilde{\Lambda}(\omega)$ gives $(\wt{\zeta}(\omega),\wt{\theta}(\omega))(A)\!=\!1$, proving \eqref{eq:mnconv}. Item \ref{thm:lambda} is now due to the right-continuity of both sides therein.

\medskip

\noindent\textbf{Proof of \ref{thm:phys}.} Let us suppose that $\Lambda_t-\Lambda_{t-}>0$ for some $t\in[0,T)$ and define $K=\inf\lb\{ z>0:\,\mu(m_{t-}(x)=0,\,x_{t-}\in[0,z])<\frac{z}{\alpha}\rb\}$.~We are going to show that $\Lambda_t-\Lambda_{t-}=K$, $\Pi$--almost surely. Using similar arguments when $\Lambda_t-\Lambda_{t-}<0$ and taking an intersection over the countable set of such discontinuities in $J$ establishes~\ref{thm:phys}.

\medskip

\noindent\emph{Upper bound.} Let us show $\Lambda_t-\Lambda_{t-}\le K$. Fix $t_n$, $\veps_n$ with $t_n,t_n+\veps_n\in([0,T]\backslash J)\cup\{0\}$ and such that $t_n\le t<t_n+\veps_n$. For $s\ge t_n$, set $\Lambda^{N,+}_s=\frac{\alpha}{N}\sum_{i=1}^N\ind_{[t_n,s]}(\sigma^{N,i})$. This is the contribution of the $X^{N,i}$--particles to the common drift $\Lambda^N$ starting from time $t_n-$. We consider the one-phase system
\begin{equation}
\begin{split}
&X_{s}^{N,i}=X_{t_n-}^{N,i}+\Big(B_s^{i}-B_{t_n}^{i}+\frac{\beta}{N}\sum_{i=1}^N{\ind_{[t_n,s]}(\tau^{N,i})}\Big)-\Lambda^{N,+}_s,\quad 1\le i\le N,\\
&\sigma^{N,i}=\inf\{s\ge 0:\,X_s^{N,i}\le 0\},\quad 1\le i\le N,\\
&\Lambda^{N,+}_s=\frac{\alpha}{N}\sum_{i=1}^N \ind_{[t_n,s]}(\sigma^{N,i}),
\end{split}
\end{equation}
where the particle $X^{N,i}$ is driven by $s\mapsto B_s^{i}-B_{t_n}^{i}+\frac{\beta}{N}\sum_{i=1}^N \mathbf{1}_{[t_n,s]}(\tau^{N,i})$
~(see Remark \ref{rmk:comparison} with $G_s^i:= B_s^{i}-B_{t_n}^{i}+\frac{\beta}{N}\sum_{i=1}^N \mathbf{1}_{[t_n,s]}(\tau^{N,i})$ ).
~We have that $((X^{N,i})_{1\le i\le N},\Lambda^{N,+})$ is a physical solution to the above system, in the sense that cascades are resolved in the manner laid out in Section \ref{sec:part} (for one-phase systems see also \cite{cuchiero_propagation_2020}, for example).

\medskip

To use the comparison principle of Remark \ref{rmk:1phasecompare}, we introduce the physical solution $((\wh{X}^{N,i})_{1\le i\le N},\wh{\Lambda}^{N})$ to the one-phase system
\begin{equation}
\begin{split}
&\wh{X}_{s}^{N,i}=X_{t_n-}^{N,i}+(B_s^{i}-B_{t_n}^{i})-\wh{\Lambda}^{N}_s,\quad 1\le i\le N,\\
&\wh{\sigma}^{N,i}=\sigma^{N,i}\wedge \inf\{s\ge t_n:\,\wh{X}_s^{N,i}\le 0\},\quad 1\le i\le N,\\
&\wh{\Lambda}^{N}_s=\frac{\alpha}{N}\sum_{i=1}^N \ind_{[t_n,s]}(\wh{\sigma}^{N,i}).
\end{split}
\end{equation}

We introduce the subprobability measures $\rho_{s-}^N:=\frac{1}{N}\sum_{i=1}^N \delta_{X_{s-}^{N,i}}\,\ind_{\{\sigma^{N,i}\ge s\}}$, $s\in[t_n,T)$, $\wh{\rho}_{s-}^N\!:=\!\frac{1}{N}\sum_{i=1}^N \delta_{\wh{X}_{s-}^{N,i}}\,\mathbf{1}_{\{\wh{\sigma}^{N,i}\ge s\}}$, $s\!\in\!(t_n,T)$, and $\wh{\rho}_{t_n-}^N\!:=\!\rho_{t_n-}^N$.~Since $\frac{\beta}{N}\sum_{i=1}^N \mathbf{1}_{[t_n,s]}(\tau^{N,i})\!\ge\! 0$, we have 
$\wh{\Lambda}^{N}_s\!\ge\! \Lambda^{N,+}_s \!\ge\! \Lambda^N_s-\Lambda^N_{t_n-}$ for all $s\!\ge\! t_n$. Following \cite[proof of Lemma~B.2]{cuchiero_propagation_2020} we find some $C>0$ such that for every $\veps>0$,
\begin{equation}\label{ineq:B2}
\begin{split}
&\;\bP\big(\forall\,0\le z\le \Lambda^N_{t_n+\veps}-\Lambda^N_{t_n-}-C\veps^{1/3}:\;\rho_{t_n-}^N([0,\alpha z+C\veps^{1/3} ])\ge z\big) \\
&\ge \bP\big(\forall\,0\le z\le \wh{\Lambda}^{N}_{t_n+\veps}-C\veps^{1/3}:\;\rho_{t_n-}^N([0,\alpha z+C\veps^{1/3} ])\ge z\big)
\ge 1-C\veps^{1/3}
\end{split}
\end{equation}
whenever $N\ge \veps^{-1/3}$. With this inequality established, we no longer need the one-phase system $((\wh{X}^{N,i})_{1\le i\le N},\wh{\Lambda}^{N})$.

\medskip

We continue by working through the proof of \cite[Theorem 6.4]{cuchiero_propagation_2020}, starting with Step 1 and noting that our crossing property (Lemma \ref{lem:cross}) implies the condition (5.3) used there. Since $t_n\in([0,T]\backslash J)\cup\{0\}$, we have
\begin{align}
\lim_{N\to \infty} \widetilde{\mu}^N_{t_n-} =  \widetilde{\mu}(m_{t_n-}(x)=0,\,x_{t_n-}\in \cdot)
\end{align}
almost surely, where the limit is taken in the space of subprobability measures endowed with the topology of weak convergence.

\medskip

Following \cite[proof of Theorem 6.4, Step 2]{cuchiero_propagation_2020}, with $\veps_n>0$ sufficiently small and $t_n+\veps_n\notin J$, we take $N\to\infty$ in \eqref{ineq:B2} to obtain
\begin{align}
\widetilde{\mu}(m_{t_n-}(x)=0,\,x_{t_n-}\in [0,\alpha z+C\veps_n^{1/3} ])\ge z(1-C\veps_n^{1/3})
\end{align}
for all $0\le z<\widetilde{\Lambda}_{t_n+\veps_n}-\widetilde{\Lambda}_{t_n-}-C\veps_n^{1/3}$, $\Pi$--almost surely. At this point, taking $n\to\infty$ and repeating the argument in \cite[proof of Theorem 6.4, Step 3]{cuchiero_propagation_2020} we obtain $\widetilde{\Lambda}_t-\widetilde{\Lambda}_{t-}\le\inf\lb\{ z>0:\,\widetilde{\mu}(m_{t-}(x)=0,\,x_{t-}\in [0,z])<\frac{z}{\alpha}\rb\}$, $\Pi$--almost surely. 

\medskip

\noindent\emph{Lower bound.} Finally, we show $L:=\Lambda_t-\Lambda_{t-}\ge K$. For the sake of contradiction, suppose $0<L< K$; we are going to see that this contradicts the right-continuity of $\Lambda$ at $t$. Define $\Lambda_{t,\veps}=\Lambda_{t+\veps}-\Lambda_t$ and note that
\begin{align}
\Lambda_{t,\veps}=\alpha(\E^\mu[m_{t+\veps}]-\E^\mu[m_t])-\beta(\E^\nu [n_{t+\veps}]-\E^\nu [n_{t}]).
\end{align}
The right-continuity of $\Lambda$ gives $\Lambda_{t,\veps}\ge -L/2$ for all $0<\veps\le \veps_0$ when $\veps_0>0$ is sufficiently small, which we shall assume. Due to the jump by $L$, we have:
\begin{equation}
\begin{split}
\E^\nu[n_{t+\veps}]-\E^\nu[n_{t}]
&=\int_{-\infty}^{-L} \theta\big(z+\sup_{t\le s \le t+\veps}(w_s-w_t-(\Lambda_{s}-\Lambda_t))\ge 0\big)\, \nu_t(\dd z)\\
&\le \int_{-\infty}^{-L} \theta\bigg(z+\sup_{t\le s \le t+\veps}(w_s-w_t)\ge -\frac{L}{2}\bigg)\, \nu_t(\dd z)\\
&\le \theta\bigg(\sup_{t\le s \le t+\veps}(w_s-w_t)\ge \frac{L}{2}\bigg)
=2\Phi\lb(-\frac{L}{2\sqrt{\veps}}\rb),
\end{split}
\end{equation}
where we have used the reflection principle for Brownian motion. Hence,
\begin{align}\label{ineq:lambdaLB}
\Lambda_{t,\veps}+2\beta\Phi\lb(-\frac{L}{2\sqrt{\veps}}\rb) \ge \alpha(\E^\mu[m_{t+\veps}]-\E^\mu[m_t]).
\end{align}
\smallskip

Using the definition of $K$ and $\Lambda_t-\Lambda_{t-}=L$, we have $\mu(m_{t}(x)=0,\,x_{t}\in [0,z])\ge\frac{z}{\alpha}$ whenever $0\le z\le K-L$. Next, we repeat the calculation from \cite[proof of Proposition 1.2]{hambly_mckean--vlasov_2018}:
\begin{equation}
\begin{split}
\E^\mu[m_{t+\veps}]-\E^\mu[m_{t}]&=\int_{0}^{\infty} \zeta\big(z+\inf_{t\le s \le t+\veps}\lb(b_s-b_t-(\Lambda_{s}-\Lambda_t)\rb)\le 0\big)\,\mu_t(\dd z)\\
&\ge \frac{1}{\alpha}\int_{0}^{K-L} \zeta\lb(b_{t+\veps}-b_t\le \Lambda_{t,\veps}-z\rb)\,\dd z\\
&= \frac{\sqrt{\veps}}{\alpha}\int_{-\Lambda_{t,\veps}/\sqrt{\veps}}^{(K-L-\Lambda_{t,\veps})/\sqrt{\veps}} \Phi(-y)\,\dd y\\
&=  \frac{\sqrt{\veps}}{\alpha} \big[x\Phi(-x)-\varphi(x)\big]\big|_{-\Lambda_{t,\veps}/\sqrt{\veps}}^{(K-L-\Lambda_{t,\veps})/\sqrt{\veps}} \\
& \ge \frac{\sqrt{\veps}}{\alpha} \big[x\Phi(-x)-\varphi(x)\big]\big|_{-\Lambda_{t,\veps}/\sqrt{\veps}}^{(K-L)/(2\sqrt{\veps})},
\end{split}
\end{equation}
where the final inequality holds when we take $\veps_0>0$ so small that $K-L-\Lambda_{t,\veps}>(K-L)/2$ for all $0<\veps\le \veps_0$.

\medskip

Note that $\varphi(x)=\varphi(-x)$ and $\Phi(-x)=1-\Phi(x)$ for all $x\in\R$. Putting this together with \eqref{ineq:lambdaLB} and simplifying yields
\begin{equation}
\begin{split}\label{ineq:Psi}
\,\varphi\lb(\frac{\Lambda_{t,\veps}}{\sqrt{\veps}}\rb)-\frac{\Lambda_{t,\veps}}{\sqrt{\veps}}\, \Phi\lb(-\frac{\Lambda_{t,\veps}}{\sqrt{\veps}}\rb)
\le \frac{2\beta}{\sqrt{\veps}}\, \Phi\lb(-\frac{L}{2\sqrt{\veps}}\rb)+\varphi\lb(\frac{K\!-\!L}{2\sqrt{\veps}}\rb)-\frac{ K\!-\!L}{2\sqrt{\veps}}\, \Phi\lb(-\frac{K\!-\!L}{2\sqrt{\veps}}\rb).
\end{split}
\end{equation}
Letting $\Psi(x):=\varphi(x)-x\Phi(-x)$ for $x\in\R$, we recognize that the left-hand side of~\eqref{ineq:Psi} can be written as $\Psi\lb(\Lambda_{t,\veps}/\sqrt{\veps}\rb)$. We note that $\Psi'(x)=-\Phi(-x)$, $\Psi$ is a decreasing convex function, $\lim_{x\to-\infty}\Psi(x)=\infty$, and $\lim_{x\to\infty}\Psi(x)=0$.~The right-hand side of~\eqref{ineq:Psi} goes to zero as $\veps\downarrow 0$; therefore, we must have $\lim_{\veps\downarrow0}\Lambda_{t,\veps}/\sqrt{\veps}=\infty$.

\medskip

We can say more thanks to inequality \eqref{ineq:Psi}. For $x>0$, recall the bounds
\begin{align}
\Phi(-x)\le \frac{\varphi(x)}{x}\quad\text{and}\quad
\Phi(-x)\le \lb( \frac{1}{x}- \frac{1}{x^3}+ \frac{1}{x^5}\rb) \varphi(x),
\end{align}
which can be obtained by (repeated) integration by parts.~Moving $\frac{\Lambda_{t,\veps}}{\sqrt{\veps}}\,\Phi\big(\!-\frac{\Lambda_{t,\veps}}{\sqrt{\veps}}\big)$ to the right-hand side of \eqref{ineq:Psi} and using these bounds for $\Phi$ yields the estimate
\begin{equation}
\begin{split}
\varphi\lb(\frac{\Lambda_{t,\veps}}{\sqrt{\veps}}\rb)
\le \frac{4\beta}{L}\,\varphi\lb(-\frac{L}{2\sqrt{\veps}}\rb)+\varphi\lb(\frac{K-L}{2\sqrt{\veps}}\rb)+\lb( 1- \frac{\veps}{(\Lambda_{t,\veps})^2}+ \frac{\veps^2}{(\Lambda_{t,\veps})^4}\rb)  \varphi\lb(\frac{\Lambda_{t,\veps}}{\sqrt{\veps}}\rb).
\end{split}
\end{equation}
Equivalently,
\begin{equation}
\begin{split}
\frac{\veps}{(\Lambda_{t,\veps})^2}\,\varphi\lb(\frac{\Lambda_{t,\veps}}{\sqrt{\veps}}\rb)
\le \frac{4\beta}{L}\,\varphi\lb(-\frac{L}{2\sqrt{\veps}}\rb)+\varphi\lb(\frac{K-L}{2\sqrt{\veps}}\rb)+ \frac{\veps^2}{(\Lambda_{t,\veps})^4}\,\varphi\lb(\frac{\Lambda_{t,\veps}}{\sqrt{\veps}}\rb).
\end{split}
\end{equation}
However, taking $\veps>0$ small enough, and since $\Lambda_{t,\veps}\to 0$, $\Lambda_{t,\veps}/\sqrt{\veps}\to\infty$ as $\veps \to 0$, we are able to make each term on the right-hand side strictly less than one third of the left-hand side.~This contradiction proves that $\Lambda$ cannot be right-continuous at $t$ if $0<\Lambda_t-\Lambda_{t-}=L< K$. \qed

\section{Two-Region Systemic Risk Model: Conclusions and Extensions}\label{sec:conc}

\begin{figure}
    \hspace{0.7cm}
    \includegraphics[height=9.5cm]{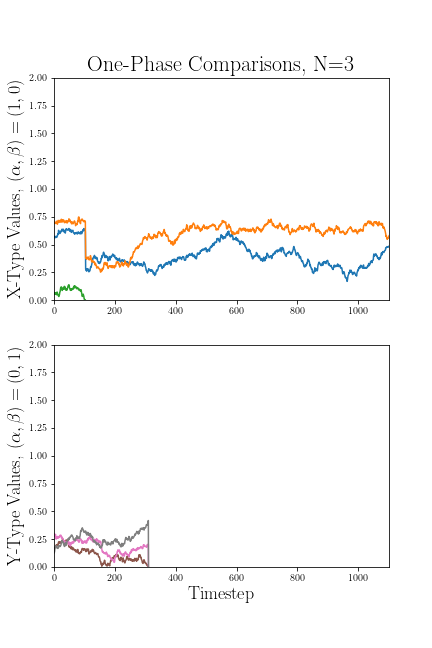}
    \includegraphics[height=9.5cm]{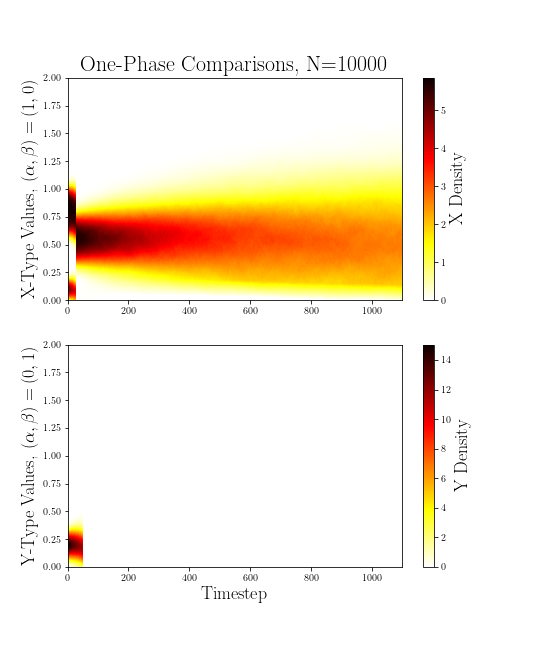}
    \caption{One-phase comparisons for \eqref{eq:tworegion} with initial densities given by \eqref{eq:exampledenisties}.}
    \label{fig:plotscomparison}
\end{figure}

We now return to the two-region systemic risk model from Section \ref{subsec:app} and the example with initial conditions given by \eqref{eq:exampledenisties}. In Figure \ref{fig:plotscomparison}, we employ the reasoning from Remark \ref{rmk:1phasecompare} to compare the firms from each region to one-phase models similar to those in \cite{nadtochiy_particle_2019} (``similar'' since the assets there are driven by geometric Brownian motions rather than Brownian motions, but this is not an issue). We picture here the same realizations of the driving Brownian motions as in Figure \ref{fig:plots}, but without the couplings between the two regions: the $\scrX$-type firms here experience an interaction term with $(\alpha,\beta)=(1,0)$, and $\scrY$-type firms experience an interaction term with $(\alpha,\beta)=(0,1)$. For the given initial conditions, we see that the $\scrY$-type firms originating in $[0.10,0.30]$ do not experience a boost from the default of the $\scrX$-type firms which start in $[0.05,0.15]$. As a result, it is not long before there is a simultaneous default of all $\scrY$-type firms (by chance, this happens somewhat later for the realization with $N=3$). Using the comparison principle of Remark \ref{rmk:comparison}, we see in general that the firms from both regions receive a more negative effect from the interaction terms in the absence of the coupling between the regions. 

\medskip

For the model considered here, defaults in one region can delay defaults in the other region, and the coupling makes each region more stable than in a one-phase model with the same initial conditions. This is due to the sign conventions in the definition of $\Lambda^N$ in \eqref{eq:tworegion} and $\Lambda$ in the limiting problem \eqref{limitprob}. In the one-phase model from \cite{nadtochiy_particle_2019}, the parameter $\alpha$ ($C$ in that work) is interpreted as an exposure parameter between firms in the same region. A systemic event, where a macroscopic proportion of firms default, is heralded by
the physical jump condition for the interaction term. In turn, the physical jump condition for the one-phase model is determined solely by the parameter $\alpha$ and the CDF of the solvent banks with assets near the default level.
For the two-phase model, a systemic event for region $\scrX$ is triggered by the same quantities but may also be averted by defaults in region $\scrY$, which depend on the CDF of $\scrY$-type firms near the default level and the parameter $\beta$. We see that $\beta$ is not only a measure of the exposure between $\scrY$-type banks, but also measures the degree to which there is a ``flight-to-quality'' from region $\scrY$ to region $\scrX$ (a similar interpretation applies vice versa for $\alpha$). The importance of the degree of interconnection within a network to the potential for systemic risk and financial contagion is well studied (see, for instance, \cite{fouque_network_2013}). Adding to this, we see with the present model that flight-to-quality between regions may provide a countering force to prevent a systemic event. 

\medskip

To separate out the two roles of $\alpha$ and $\beta$, the model may be extended by replacing the common interaction term with 
\begin{align}\label{eq:decouple}
\Lambda^{\scrX}_t=\alpha^{\scrX} \bP(\sigma\le t)-\beta^{\scrX} \bP(\tau\le t)\quad\text{and}\quad\Lambda^{\scrY}_t=\alpha^{\scrY} \bP(\tau\le t)-\beta^{\scrY} \bP(\sigma\le t)
\end{align}
for $\scrX$- and $\scrY$-type firms, respectively. Here, $\alpha^{\scrX}$ and $\alpha^{\scrY}$ measure exposure within the respective regions. Meanwhile, $\beta^{\scrX}$ measures the flight-to-quality effect from region $\scrY$ to $\scrX$ (and vice versa for $\beta^{\scrY}$). Recalling Remark \ref{rmk:alphabeta}, we may also consider $\alpha^{\scrX}>0$ and $\beta^\scrX<0$, which renders $t\mapsto \Lambda^{\scrX}$ monotone.
A negative $\beta^{\scrX}$ causes defaults in  region $\scrY$ to propagate to region $\scrX$. Correlation between regions in times of financial crisis is well known from empirical economics literature.
A need for models with cross-exciting jumps which capture this behaviour has been discussed in works such as \cite{ait-sahalia_modeling_2015,ait-sahalia_portfolio_2016} where a generalization of the Hawkes process is used.
In contrast to existing models where the jumps are driven by an exogenous intensity process, the jumps in the present model are endogenous as they are triggered by the hitting times which correspond to defaults.

\medskip

The systemic risk model considered in this work, despite its relative simplicity, gives rise to interesting behaviors which can be seen in real-world financial networks. Going from one region to two regions presented difficulties due to the non-monotone interaction term $\Lambda$. The model can be generalized from two to an arbitrary number of regions without difficulty: just as with equation \eqref{eq:decouple}, each region can be given an interaction term which is a linear combination of the CDFs of the hitting times from all regions. Some bookkeeping needs to be done as well to resolve possible cascades in a physical manner, as in Section \ref{sec:part}.
~The driving processes for the individual firms can be given by any c\`adl\`ag paths, as mentioned in Remark \ref{rmk:comparison} as well as \cite{nadtochiy_mean_2019}. In the financial context, the case of Brownian motions with a common noise, studied in \cite{hambly_spde_2019} for a one-phase model, is a natural next step. 
In another direction, just as the one-phase systemic risk model poses an optimal bailout problem for a regulator seeking to prevent a systemic event, such as in \cite{cuchiero_optimal_2021}, the bailout problem with one regulator in each region poses a game theoretic problem of potential interest.


\bigskip\bigskip

\printbibliography

\bigskip\bigskip

\end{document}